\documentclass[12pt]{article}

\usepackage[dvips]{graphicx}

\usepackage{amssymb,amsfonts,amsmath,amsthm}
\usepackage{color}

\usepackage{url}

\newtheorem{theorem}{Theorem}[section]
\newtheorem{lemma}[theorem]{Lemma}
\newtheorem{proposition}[theorem]{Proposition}
\newtheorem{corollary}[theorem]{Corollary}

\theoremstyle{definition}

\newtheorem{example}[theorem]{Example}

\theoremstyle{remark}
\newtheorem{remark}[theorem]{Remark}

\numberwithin{equation}{section}

\newcommand{\h}{{\mathcal H}}
\newcommand{\mn}{\mathbb N}

\newcommand{\mr}{\mathbb R}

\def\range{{\mathcal R}}

\def\LtRd{L^2\left(\mathbb{R}^d\right)}

\def\h{{\mathcal H}}
\def\bp{\noindent{\bf Proof: \ }}
\def\ep{\noindent{$\Box$}}

\def\<{\langle}
\def\>{\rangle}

\newcommand{\seq[1]}{(#1_n)}

\def\newin {\,\kern-0.4em\in\kern-0.15em}
\def\newsubset {\kern-0.2em\subset\kern-0.2em}

\def\v{\vspace{.1in}}

\title{Representation of the inverse of a frame multiplier}

\author{Peter Balazs and Diana T. Stoeva \vspace{.1cm} \\ 
  Acoustics Research Institute, Austrian Academy of Sciences,\\
Wohllebengasse 12-14, Vienna A-1040, Austria \\
}

\begin{document}

\maketitle

\begin{abstract}
  Certain mathematical objects appear in a lot of scientific disciplines, like physics, signal processing and, naturally, mathematics. In a general setting they can be described as frame multipliers, consisting of analysis, multiplication by a fixed sequence (called the symbol), and synthesis. They are not only interesting mathematical objects, but also  important for applications, for example for the realization of time-varying filters. In this paper we show a surprising result about the inverse of such operators, if existing, as well as new results about a core concept of frame theory, dual frames. 
We show that for semi-normalized symbols, the inverse of any invertible frame multiplier can always be represented 
as a frame multiplier with the reciprocal symbol and dual frames of the given ones. 
Furthermore, one of those dual frames is uniquely determined and the other one can be arbitrarily chosen. 
We investigate sufficient conditions for the special case, when both dual frames can be chosen to be the canonical duals. 
In connection to the above, we show that the set of dual frames determines a frame uniquely. 
Furthermore, for a given frame, the union of all coefficients of its dual frames is dense in $\ell^2$. 
We also investigate invertible Gabor multipliers. 
 Finally, we give a numerical example for the invertibility of multipliers in the Gabor case.
 \end{abstract}



\section{Introduction, Notation, and Motivation}

In many scientific disciplines, certain objects play an important role. 
Those systems are described by an analysis procedure followed by a multiplication, followed by a synthesis. 
Those operators are of utmost importance in 
\begin{itemize}
\item Mathematics, where they are used for the diagonalization of operators \cite{schatt1};
\item Physics, where they are a link between classical and quantum mechanics, so called quantization operators \cite{aliant1};
\item Signal processing, where they are a particular way to implement time-variant filters \cite{hlawatgabfilt1}; 
\item Acoustics, where those time-frequency filters are used in several fields, for example in Computational Auditory Scene Analysis \cite{wanbro06}. 
\end{itemize}

In this paper we show a surprising result about the shape of the inverse of such operators, if existing. 
This also leads us to new results concerning dual frames, a concept at the core of frame theory. 

To be able to describe those operators in a general setting, as an extension of Gabor multipliers \cite{feinow1},
multipliers for general Bessel sequences were introduced by one of the authors \cite{xxlmult1}.
Further, multipliers for general sequences were investigated in \cite{stoevxxl09,balsto09new,uncconv2011,bstable09}.  
These are operators defined by
\begin{equation} \label{multdef}
 M_{m,\Phi, \Psi}h=\sum_{n=1}^\infty m_n \<h,\psi_n\>\phi_n ,
\end{equation} 
for given 
sequences $\Phi=(\phi_n)$ and $\Psi=(\psi_n)$ with elements from a Hilbert space $\h$, and 
a given complex scalar sequence $m=(m_n)$ called the {\it symbol}.
Such operators are also investigated for continuous transforms - in a general \cite{xxlbayasg11} (continuous frame multipliers), wavelet \cite{rochberg90} (Calderon-Toeplitz operators) and short-time Fourier setting \cite{cogr05} (localization operators).
Here we stick to the discrete version.
 Multipliers  are interesting not only from a theoretical point of view, but also for applications. 
They are applied for example in psychoacoustical modelling \cite{xxllabmask1}
and denoising \cite{majxxl10}. 
Multipliers are a particular way to implement time-variant filters
 \cite{hlawatgabfilt1}.
Therefore, for some applications it is important to find the inverse of a multiplier if it exists.
 The paper \cite{balsto09new} is devoted to invertibility of multipliers, necessary conditions for invertibility, sufficient conditions, and representation for the inverse via Neumann series. 

  In the present paper our attention is on how to express the inverse of an invertible frame multiplier as a multiplier with the reciprocal symbol and dual frames of the given ones.
 We show a result for all frames, 
 namely, the inverse of any invertible frame multiplier with a semi-normalized symbol 
can always be represented as a multiplier with the reciprocal symbol and dual frames of the given ones, 
 where one of these dual frames is uniquely determined and the other one can be arbitrarily chosen:
\begin{theorem} \label{ff} 
Let $\Phi$ and $\Psi$ be frames for $\h$, and let 
the symbol  
$m$ 
satisfy $0<\inf_n |m_n|\leq \sup_n|m_n|<\infty$. 
Assume that $M_{m, \Phi, \Psi}$ is invertible. 
Then 
\begin{itemize}
\item
there exists a unique dual frame $\Phi^\dagger$ of $\Phi$, 
so that for any dual frame $\Psi^d$ of $\Psi$ we have
\begin{equation}\label{minv1}
M_{m, \Phi, \Psi}^{-1}  = M_{1/m, \Psi^d, \Phi^\dagger};
\end{equation} 
\item
there exists a unique dual frame $\Psi^\dagger$ of $\Psi$, 
so that for any dual frame $\Phi^d$  of $\Phi$  we have
\begin{equation}\label{minv2}
M_{m, \Phi, \Psi}^{-1}  = M_{1/m, \Psi^\dagger, \Phi^d}.
\end{equation} 
\end{itemize}

\end{theorem}

The investigation of this topic led us to surprising results about dual frames. 
We show that a frame is uniquely defined by the set of its dual frames. Furthermore, for a given frame, the union of all coefficients of its dual frames is dense in $\ell^2$:
   \begin{theorem}\label{lemrange} 
  Let $\Phi$ be a frame for $\h$. 
  Then the following statements hold.
  \begin{itemize}
   \item[{\rm (i)}] 
   The closure of the union of all sets $\range(U_{\Phi^d})$, where $\Phi^d$  runs \,
   through all dual frames of $\Phi$, is $\ell^2$.  
   \item[{\rm (ii)}] Let  $\Psi$ be a frame for $\h$. If 
 every dual frame $\Phi^d$ of $\Phi$ is a dual frame of $\Psi$, then $\Psi=\Phi$.
  \end{itemize}
  \end{theorem}

\begin{figure}[ht] 
\includegraphics[width=1\textwidth]{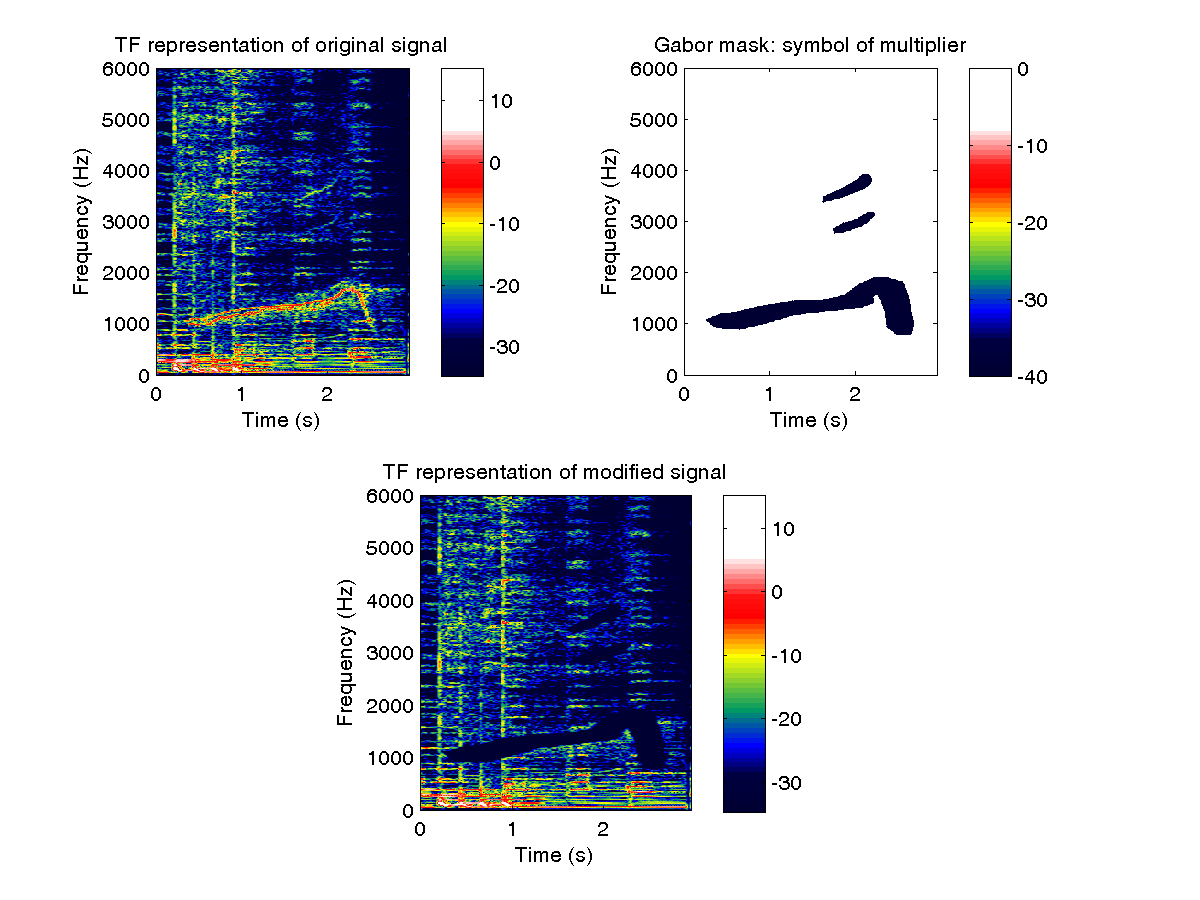} 
\caption{\label{fig:examp1} {\em An illustrative example to visualize a multiplier.} 
(TOP LEFT) The time-frequency representation of 
 the music signal $f$. (TOP RIGHT) The symbol $m$, found by a (manual) estimation of the time-frequency region of the singer's voice. 
 (BOTTOM)  Time-frequency representation of $M_{m,\widetilde \Psi,\Psi}f$.} 
\end{figure}

In Figure \ref{fig:examp1} we show a visualization of a multiplier $M_{m,\Phi,\Psi}$ in the time-frequency plane, 
 which will again become interesting  in the last section of the paper. The visualization is done using algorithms in the LTFAT toolbox \cite{soendxxl10}. We consider a music signal $f$ and the action of a multiplier $M_{m,\Phi,\Psi}$ on $f$. 
For $f$ we use a $2$ seconds long excerpt of the ``Jump'' from Van Halen. 
For a time-frequency representation of the musical signal $f$ (TOP LEFT) we use a 'painless' Gabor frame $\Psi$ (a $80$ ms Hanning window with $12,5$\% overlap). 
By a manual estimation, we determine the symbol $m$ that should describe the time-frequency region of the singer's voice. This region is then multiplied by $0.01$, the rest by $1$ (TOP RIGHT). 
Finally, we show a time-frequency representation of the modified signal (BOTTOM).

For implementations and scripts producing Figures 1 and 2, see 
\url{http://www.kfs.oeaw.ac.at/RepresentationInverseMultiplier.}

\subsubsection*{Motivation} 

In \cite{xxlmult1} it is proved that, if $m$ is semi-normalized, then 
a Riesz multiplier  $M_{m,\Phi,\Psi}$ is automatically invertible  and
\begin{equation}\label{invmrb}
M_{m,\Phi,\Psi}^{-1}=M_{1/m,\widetilde{\Psi},\widetilde{\Phi}},
\end{equation} 
where $\widetilde{\Phi}$ and  $\widetilde{\Psi}$ denote the canonical duals of $\Phi$ and $\Psi$, respectively. 

The result on Riesz multipliers has opened the following questions: 
\begin{quote} 
 [{$\bf Q1$}] {\em 
Are there other invertible frame multipliers $M_{m,\Phi,\Psi}$ whose inverses can be represented using the inverted symbol $1/m$ and appropriate dual frames of $\Phi$ and $\Psi$?\\
 } 
\end{quote}

\begin{quote} 
[{$\bf Q2$}] {\em
 Are there other invertible frame multipliers $M_{m,\Phi,\Psi}$ whose inverses can be written as $M_{1/m,\widetilde{\Psi},\widetilde{\Phi}}$ using the canonical duals 
 }? 
\end{quote}
 
 \v The paper is devoted to these two questions. 
First note that every bounded (resp. bounded surjective) operator can be written as a Bessel (resp. frame) multiplier.  Thus, the inverse of every invertible multiplier can be written as a frame multiplier. 
The aim of the present paper is to represent the inverse of an invertible frame multiplier as described in  [$Q1$].

 We give an affirmative answer to Question [$Q1$]. 
 We show  in  Theorem \ref{ff} that the inverse of every invertible frame multiplier with semi-normalized symbol can be
represented as a multiplier with the reciprocal symbol and dual frames of the given ones. One
of the dual frames is uniquely determined, while the other one can be arbitrarily
chosen.

 We also give an affirmative answer to Question [$Q2$]. We determine frame multipliers $M_{m,\Phi,\Psi}$ (not necessarily being Riesz multipliers) which are invertible and their inverses can be written as $M_{1/m,\widetilde{\Psi},\widetilde{\Phi}}$.

 The last part of this paper is devoted to Gabor multipliers.  
We determine equivalent conditions for an invertible operator on $\LtRd$ 
(and its inverse)
 to be represented as a Gabor frame multiplier with a constant symbol.

\subsubsection*{Notation and definitions}

Throughout the paper, $\h$ denotes a separable Hilbert space, $\Phi=(\phi_n)_{n=1}^\infty$ and $\Psi=(\psi_n)_{n=1}^\infty$ are sequences with elements from $\h$. The sequence $(e_n)_{n=1}^\infty$ denotes an orthonormal basis of $\h$ and $(\delta_n)_{n=1}^\infty$ denotes the canonical basis of $\ell^2$.
When the index set is omitted, $\mn$ should be understood as the index set.
The letter $m$ is used to denote a complex valued scalar sequence $(m_n)$. 
Furthermore, $\overline{m}=(\overline{m}_n)$ and $1/m=(1/m_n)$. 
The sequence $m$ is called {\it semi-normalized} if $0<\inf_n |m_n|\leq \sup_n|m_n|<\infty$. 
 For $m\in\ell^\infty$, we will use the operator $\mathcal{M}_m:\ell^2\to\ell^2$ given by $\mathcal{M}_m \seq[c]=(m_nc_n)$, which is bounded with $\| \mathcal{M}_m\|=\|m\|_{\ell^\infty}$. 
An operator $M :\h\to \h$ is called {\it invertible}  if it is a bounded bijection  from $\h$ onto $\h$. The identity operator on $\h$ is denoted by ${\rm Id}_\h$.

Recall that $\Phi$ is called a {\it frame for $\h$ with bounds $A_\Phi, B_\Phi$} if $0<A_\Phi\leq B_\Phi < \infty$ 
and $A_\Phi\|h\|^2\leq \sum_{n=1}^\infty |\<h,\phi_n\>|^2 \leq B_\Phi\|h\|^2$ for every $h\in\h$. 
For a given frame $\Phi$ for $\h$, the analysis operator is denoted by $U_\Phi$, the synthesis operator by $T_\Phi$, the frame operator by $S_\Phi$, a dual frame of $\Phi$ by $\Phi^d=(\phi^d_n)$, and the canonical dual by $\widetilde{\Phi}=(\widetilde{\phi}_n)$. For the definition of all these frame-related concepts, as well as for the definition of a Bessel sequence and a Riesz basis, we refer to \cite{ole1}.
Recall that two frames $\Phi$ and $\Psi$ for $\h$ are called {\em equivalent} if there exists an invertible operator $G:\h\to\h$ so that $\psi_n=G\phi_n$ for all $n\in\mn$. When $\Phi$ is a frame for $\h$, then a dual frame $\Phi^d$ of $\Phi$ is equivalent to $\Phi$ if and only if $\Phi^d=\widetilde{\Phi}$ \cite[Sect. 1.2]{HL}.

For given $m$, $\Phi$, and $\Psi$, the operator $M_{m,\Phi,\Psi}$ given by Equation (\ref{multdef})
is called a {\it multiplier}. The operator $M_{m,\Phi,\Psi}$ is called {\it unconditionally convergent} if the series in Equation (\ref{multdef})  converges unconditionally for every $h\in\h$. When $\Phi$ and $\Psi$ are Bessel sequences, frames, Riesz bases for $\h$, then $M_{m,\Phi,\Psi}$ will be called a {\it Bessel multiplier, frame multiplier, Riesz multiplier,} respectively. 
When $m\in\ell^\infty$, then a Bessel multiplier is a well defined operator from $\h$ into $\h$ \cite{xxlmult1}.

For $\omega \in \mr^d$ and $\tau\in\mr^d$, recall the modulation operator $E_\omega:\LtRd\to\LtRd$ 
and the translation operator $T_\tau:\LtRd\to \LtRd$ given by
$$ (E_\omega f) (x)=e^{2\pi i \omega x} f(x), \ (T_\tau f) (x)= f(x-\tau).$$ 
The symbol $\Lambda=\{ (\omega, \tau)\}$  denotes a (possibly irregular) lattice in $\mr^{2d}$. 
For $\lambda=(\omega,\tau)\in\Lambda$, the operator  $E_\omega T_\tau$ is denoted by $\pi(\lambda)$. 
For a given $g \in \LtRd$, a sequence  of the form  $(g_\lambda)_{\lambda\in\Lambda} = (\pi(\lambda)g)_{\lambda\in\Lambda}$ is called a {\it Gabor system}. When a Gabor system is a frame, it is called a {\it Gabor frame}. 
Recall that the canonical dual of a Gabor frame  $(g_\lambda)_{\lambda\in\Lambda}$ is the Gabor frame
$(\widetilde{g}_\lambda)_{\lambda\in\Lambda}=(\pi(\lambda)\widetilde{g})_{\lambda\in\Lambda}$, where $\widetilde{g}=S_{(\pi(\lambda)g)_{\lambda\in\Lambda}}^{-1}g$.
When $\Phi$ and $\Psi$ are Gabor systems (resp. Gabor frames), then $M_{m,\Phi,\Psi}$ is called a {\it Gabor multiplier} (resp. {\it Gabor frame multiplier}).

\section{The set of dual frames} \label{df}

In order to prove Theorem \ref{ff} and Proposition \ref{ffff}, we need Theorem \ref{lemrange}, stated on page \pageref{lemrange}.
This is a result which is of independent interest for frame theory, showing new properties of the set of dual frames.

  {\noindent{\bf Proof of Theorem \ref{lemrange}: \ }}
 (i)  Let the sequence $c=(c_n)\in\ell^2$ fulfill 
   $c\perp  \range(U_{\Phi^d})$ for every dual frame $\Phi^d$ of $\Phi$. 
  Then 
  \begin{equation} \label{synth}
  T_{\Phi^d} c =0, \ \forall \mbox{ dual frame } \Phi^d  \mbox{ of } \Phi.
  \end{equation} 
  The dual frames of $\Phi$ are precisely the sequences 
  $$(\widetilde{\phi}_n + h_n - \sum_{j=1}^\infty \<\widetilde{\phi}_n, \phi_j\>h_j)_{n=1}^\infty, $$
  where $(h_n)_{n=1}^\infty$ is a Bessel sequence in $\h$ 
  (see, e.g., \cite[Theorem 5.6.5]{ole1}).
  Therefore, 
  \begin{equation*}
\sum_{n=1}^\infty c_n \left(\widetilde{\phi}_n + h_n - \sum_{j=1}^\infty \<\widetilde{\phi}_n, \phi_j\>h_j\right) =0
  \end{equation*}
  for every  Bessel sequence $\{h_n\}_{n=1}^\infty$ in $\h$.
  By Equation (\ref{synth}) we have $ T_{\widetilde{\Phi}} c =0$, which implies that 
   \begin{equation} \label{hbessel}
\sum_{n=1}^\infty c_n \left(h_n - \sum_{j=1}^\infty \<\widetilde{\phi}_n, \phi_j\>h_j\right) =0
  \end{equation}
  for every Bessel sequence $\{h_n\}_{n=1}^\infty$ in $\h$. 
  Using Equation (\ref{hbessel}) with the Bessel sequence $(h_n)_{n=1}^\infty=(e_1, 0,0,0,\ldots) $, 
we obtain
   \begin{equation*}
 c_1 e_1 - \sum_{n=1}^\infty c_n \<\widetilde{\phi}_n, \phi_1\>e_1=0,
  \end{equation*}
which implies that
 $c_1=0$.
\sloppy
  In a similar way, using Equation (\ref{hbessel}) with the Bessel sequence $(h_n)_{n=1}^\infty=(0, \ldots, 0, e_j, 0,0,0,\ldots) $, where $e_j$ stands at the $j$-th position, we obtain $c_j=0$ for every $j\geq 2$.
 Therefore, $c=(0)$, which completes the proof.
 
 (ii) 
Assume that all dual frames $\Phi^d$ of $\Phi$ are dual frames of $\Psi$. Then $T_\Phi U_{\Phi^d}= Id_\h=T_\Psi U_{\Phi^d}$, which by (i) implies that $T_\Phi =T_\Psi$ and hence, $\Phi=\Psi$. 
    \ep

\v 
By the above result, different frames have different sets of dual frames;
if two frames $\Phi$ and $\Psi$ for $\h$ have the same sets of dual frames, then $\Phi=\Psi$. 
In particular, two different frames cannot have sets of dual frames which are included into one another.

\section{Inversion of Multipliers by Inverted Symbol [Q1] and Dual Frames} \label{secq1}

Here we give an affirmative answer to Question [$Q1$]. 
The result about the inverses of invertible frame multipliers with semi-normalized symbols is stated in Theorem \ref{ff}. 
In addition, we show the following:

\begin{proposition} \label{ffff} 
For the assumptions in Theorem \ref{ff}, we have the additional properties:
\begin{itemize}
\item
If $F=(f_n)$ is a Bessel sequence in $\h$ which fulfills 
$M_{m, \Phi, \Psi}^{-1}  = M_{1/m, \Psi^\dagger, F}$
(resp. $M_{m, \Phi, \Psi}^{-1}  = M_{1/m, F, \Phi^\dagger}$),
then $F$ must be a dual frame of $\Phi$ (resp. $\Psi$).
\item
 $\Psi^\dagger$ is the only Bessel sequence in $\h$ which satisfies 
$M_{m, \Phi, \Psi}^{-1}  = M_{1/m, \Psi^\dagger, \Phi^d}$  
for all dual frames $\Phi^d$ of $\Phi$. 
\item 
 $\Phi^\dagger$ is the only Bessel sequence in $\h$ which satisfies 
 $M_{m, \Phi, \Psi}^{-1}  = M_{1/m, \Psi^d, \Phi^\dagger}$  
for all dual frames  $\Psi^d$ of $\Psi$. 
\end{itemize}
\end{proposition}
\v
{\noindent{\bf Proof of Theorem \ref{ff} and Proposition \ref{ffff}: \ }}
Denote $M:= M_{m, \Phi, \Psi}$. First observe that the sequence $(M^{-1} (m_n \phi_n))$ is a dual frame of $\Psi$.  Denote it by $\Psi^\dagger$. 
Therefore, $M^{-1}T_\Phi \delta_n=T_{\Psi^\dagger}\mathcal{M}_{1/m} \delta_n$, $n\in\mn$.
Now the boundedness of the operators imply that
$M^{-1} T_\Phi  =T_{\Psi^\dagger} \mathcal{M}_{1/m}$ on $\ell^2$.
Using any dual frame $\Phi^d$ of $\Phi$ we get
$ M^{-1}  =  T_{\Psi^\dagger}  \mathcal{M}_{1/m} U_{{\Phi^d}}  = M_{1/m, \Psi^\dagger, {\Phi^d}}
 \ \mbox{on} \ \h.$
 
In a similar way as above, it follows that the sequence $((M^{-1})^*(\overline{m_n}\psi_n))$ is a dual frame of $\Phi$ (denoted by $\Phi^\dagger$) and 
hence, 
\begin{equation} \label{phi1}
(M^{-1})^* T_\Psi  = T_{\Phi^\dagger} \mathcal{M}_{1/\overline{m}} \mbox{ on } \ell^2.
\end{equation}
Therefore, 
$M^{-1} =T_{{\Psi^d}} \mathcal{M}_{1/m} U_{\Phi^\dagger}=M_{1/m, {\Psi^d}, \Phi^\dagger}.$

Now assume that 
   $F=(f_n)$ is a Bessel sequence in $\h$ which satisfies 
$M_{m, \Phi, \Psi}^{-1}  = M_{1/m, F,\Phi^\dagger} $. 
By Equation (\ref{phi1}), it follows that 
$T_\Psi U_F= M^* T_{\Phi^\dagger} \mathcal{M}_{1/\overline{m}} U_F 
=M^* (M^{-1})^*
={\rm Id}_\h,$ which implies that $F$ is a dual frame of $\Psi$.
 In a similar way, every Bessel sequence $F$ in $\h$ which satisfies 
$M_{m, \Phi, \Psi}^{-1}  = M_{1/m, \Psi^\dagger, F} $ must be a dual frame of $\Phi$. 

On the other hand, assume that $F$ is a Bessel sequence in $\h$ which satisfies 
 $M_{1/m, F, \Phi^d}= M_{1/m, \Psi^\dagger, \Phi^d}$  for all dual frames $\Phi^d$ of $\Phi$.
Then $T_F  \mathcal{M}_{1/m} U_{{\Phi^d}}= T_{ \Psi^\dagger}  \mathcal{M}_{1/m} U_{{\Phi^d}}$  for all dual frames $\Phi^d$ of $\Phi$, which by Theorem \ref{lemrange}(i) implies that 
$T_F  \mathcal{M}_{1/m} = T_{ \Psi^\dagger}  \mathcal{M}_{1/m} $.
Since $m$ is semi-normalized (so, $\mathcal{M}_{1/m}$ is invertible on $\ell^2$), 
it follows that $T_F   = T_{ \Psi^\dagger}   $ and hence, $F=\Psi^\dagger$.

 The statement for $\Phi^\dagger$ follows in a similar way.
\ep

\v 
Concerning Theorem \ref{ff}, it is natural to ask whether the frame $\Psi^\dagger$ (resp $\Phi^\dagger$) is the canonical dual of $\Psi$ (resp. $\Phi$). Observe that in this context we have $\Psi^\dagger=\widetilde{\Psi}$ (resp. $\Phi^\dagger=\widetilde{\Phi}$ ) if and only if 
  $\Psi$ is equivalent
  to $(m_n\phi_n)$ (resp. $\Phi$ is equivalent to $(\overline{m}_n\psi_n)$ ).

\v
 Note that Equations (\ref{minv1}) and (\ref{minv2})
are {\em not} constructive approaches 
leading to an implementation for the inversion of $M$. For the dual frame $\Psi^\dagger$ (resp. $\Phi^\dagger$) we already had to apply $M^{-1}$. For more constructive approaches to the inversion of multipliers see the next sections and \cite{balsto09new}. 

\v 
A sub-result of Theorem \ref{ff}, 
the representation of the inverse for the particular case of finite-dimensional spaces and $\Psi = \Phi^d$, has been independently found in the context of frame diagonalization of matrices \cite{Futamura20123201}. 

\v

\begin{remark} 
In \cite{uncconv2011} the following conjecture is formulated:
For an unconditionally convergent multiplier $M_{m,\Phi,\Psi}$, 
 there always exist sequences
 $(c_n)$ and $(d_n)$ so that $M_{m,\Phi,\Psi}=M_{(1),(c_n\phi_n), (d_n\psi_n)}$ 
and the sequences $(c_n\phi_n)$, $(d_n\psi_n)$ are Bessel sequences.

If this conjecture is true,  then any invertible, unconditionally convergent multiplier $M_{m,\Phi,\Psi}$ can be rewritten as $M_{m,\Phi,\Psi}=M_{(1),(c_n\phi_n), (d_n\psi_n)}$ 
where the sequences $(c_n\phi_n)$, $(d_n\psi_n)$ are frames for $\h$,
 and thus, by Theorem \ref{ff}, $M_{m,\Phi,\Psi}^{-1}$ can be written as 
$M_{(1), (d_n\psi_n)^\dagger, (c_n\phi_n)^d}$ and 
 $M_{(1), (d_n\psi_n)^d,(c_n\phi_n)^\dagger}$.
\end{remark}

\section{Inversion of Multipliers Using the Canonical Duals [Q2]} \label{secq2}

The following example shows cases where Question [$Q2$] is answered affirmatively. 

\begin{example} \label{overcfr1} 
Every frame $\Phi$ for $\h$ fulfills $M_{(1),\Phi,\Phi}^{-1}=M_{(1),\widetilde{\Phi},\widetilde{\Phi}}$. 
 \end{example}

Example \ref{overcfr3} shows a case when $M_{m,\Phi,\Psi}$ is invertible but the inverse is not equal to $M_{1/m,\widetilde{\Psi},\widetilde{\Phi}}$.

\begin{example} \label{overcfr3}
 Let $\Phi=(e_1,e_1,e_1,e_2,e_2,e_2,e_3,e_3,e_3,\ldots)$ and 
 $\Psi=(e_1, e_1, -e_1, e_2, e_2, -e_2, e_3, e_3, -e_3,\ldots)$. 
  Then 
  $M_{(1),\widetilde{\Psi},\widetilde{\Phi}}=\frac{1}{9} \,{\rm Id}_\h \neq M_{(1),\Phi,\Psi}^{-1}={\rm Id}_\h$.
\end{example}

The next proposition determines a class of  multipliers which are invertible and whose inverses can be written as in Equation (\ref{invmrb}). 
While in Theorem \ref{ff} 
it is assumed that the frame multiplier is invertible,
in Proposition \ref{mi2} 
we investigate the invertibility of frame multipliers 
- we give sufficient conditions for invertibility and sufficient conditions for non-invertibility.
For the rest of the section the letter $c$ means a non-zero constant.

\begin{proposition} \label{mi2}
Let $\Phi$ and $\Psi$ be frames for $\h$ and $\seq[m]=(c)$. 
Then the following assertions hold.
\begin{itemize}
\item[{\rm (i)}] If $\range(U_\Phi)\subseteq \range(U_\Psi)$, then 
$M_{(1/c),\widetilde{\Psi},\widetilde{\Phi}}$  is a bounded right inverse of $M_{(c),\Phi,\Psi}$.

\item[{\rm (ii)}]  If $\range(U_\Psi)\subseteq \range(U_\Phi)$, then 
$M_{(1/c),\widetilde{\Psi},\widetilde{\Phi}}$ is a bounded left inverse of $M_{(c),\Phi,\Psi}$.

\item[{\rm (iii)}]  If $\range(U_\Phi)= \range(U_\Psi)$, then $M_{(c),\Phi,\Psi}$ is invertible and 
$M_{(c),\Phi,\Psi}^{-1}=M_{(\frac{1}{c}), \widetilde{\Psi},\widetilde{\Phi}}$.

\item[{\rm (iv)}]  If $\range(U_\Phi)\subsetneq \range(U_\Psi)$, then $M_{(c),\Phi,\Psi}$ is not invertible.

\item[{\rm (v)}]  If $\range(U_\Psi)\subsetneq \range(U_\Phi)$, then $M_{(c),\Phi,\Psi}$ is not invertible.
\end{itemize}
\end{proposition}

\bp (i) Assume that $\range(U_\Phi)\subseteq \range(U_\Psi)$. For every $h\in\h$, the element $U_\Phi S_\Phi^{-1} h$ 
can be written as $U_\Psi g^h$ for some $g^h\in\h$ and 
$$M_{(c),\Phi,\Psi} M_{(1/c),\widetilde{\Psi},\widetilde{\Phi}} h = 
T_\Phi U_\Psi S_\Psi^{-1} T_\Psi U_\Phi S_\Phi^{-1} h=
T_\Phi U_\Psi g^h= h.
$$

(ii) can be proved in a similar way as (i).

(iii) follows from (i) and (ii).

(iv) Assume that $\range(U_\Phi)\subset \range(U_\Psi)$ with $\range(U_\Phi)\neq \range(U_\Psi)$. 
By (i), the operator $M_{(1/c),\widetilde{\Psi},\widetilde{\Phi}}$  is a bounded right inverse of $M_{(c),\Phi,\Psi}$. We will prove that 
$M_{(1/c),\widetilde{\Psi},\widetilde{\Phi}}$  is not a left inverse of $M_{(c),\Phi,\Psi}$, which will imply that $M_{(c),\Phi,\Psi}$ can not be invertible. 
Consider an arbitrary element $g\in \range(U_\Psi)\setminus \range(U_\Phi)$ and write $g=U_\Psi h$ for some $h\in\h$.
 Since $\ell^2=\range(U_\Phi)\oplus \ker(T_\Phi)$, we can also write $g= U_\Phi f + d$ for some $f\in\h$ and some $d\in \ker(T_\Phi)$. 
Then
\begin{eqnarray*}
 M_{(1/c),\widetilde{\Psi},\widetilde{\Phi}} M_{(c),\Phi,\Psi} h &=& S_\Psi^{-1} T_\Psi U_\Phi S_\Phi^{-1} T_\Phi U_\Psi h \\
& =&  S_\Psi^{-1} T_\Psi U_\Phi S_\Phi^{-1} T_\Phi (U_\Phi f + d)  \\
  &=&  S_\Psi^{-1} T_\Psi (U_\Psi h - d) = h -  S_\Psi^{-1} T_\Psi d.
  \end{eqnarray*}
Observe that $S_\Psi^{-1} T_\Psi d \neq 0$, which 
implies that $M_{(1/c),\widetilde{\Psi},\widetilde{\Phi}}$ is not a left inverse of $M_{(c),\Phi,\Psi}$.
 
 (v) Assume that $\range(U_\Psi)\subset \range(U_\Phi)$ with $\range(U_\Psi)\neq \range(U_\Phi)$. 
 By (i), $M_{(1/c),\widetilde{\Psi},\widetilde{\Phi}}$  is a bounded left inverse of $M_{(c),\Phi,\Psi}$.
 In a similar way as in (iv), one can prove that $M_{(1/c),\widetilde{\Psi},\widetilde{\Phi}}$  is not a right inverse of $M_{(c),\Phi,\Psi}$, which implies that $M_{(c),\Phi,\Psi}$ can not be invertible. 
\ep

\vspace{.1in}
Concerning the statements in Proposition \ref{mi2}, note that if none of $\range(U_\Phi)$ and $\range(U_\Psi)$  is a subset of the other one, then both invertibility and non-invertibility of $M_{(c),\Phi,\Psi}$ are possible. For a case of invertibility, consider the frame multiplier $M_{(1),\Phi,\Psi}$, where $\Phi=(e_1,e_1,e_2,e_2,e_3,e_3,e_4,e_4,\ldots)$ and 
$\Psi = (\frac{1}{2}e_1, \frac{1}{2}e_1, \frac{1}{2}e_2, \frac{1}{2}e_2, \frac{1}{3}e_3, \frac{2}{3}e_3,\frac{1}{4}e_4,\frac{3}{4}e_4,\ldots )$, and thus, $M_{(1),\Phi,\Psi}$ is the Identity operator on $\h$. For a case of non-invertibility, consider the frame multiplier $M_{(1),\Phi,\Psi}$, where
  $\Phi=(e_1,e_1,e_2,e_2,e_3,e_3,e_4,e_4,\ldots)$ and $\Psi=(e_1,e_1,e_2,e_3,e_4,\ldots)$.

\begin{remark} \label{remequiv}
 Let $\Phi$ and $\Psi$ be frames for $\h$. 
The condition $\range(U_\Phi)= \range(U_\Psi)$ corresponds to $\Phi$ and $\Psi$ being equivalent frames \cite[Corollary 4.5]{Casaz1}. 
The condition $\range(U_\Phi)\subseteq \range(U_\Psi)$ is identical to $\Psi$ being partial equivalent to $\Phi$, i.e. 
to the existence of a bounded operator $Q:\h\to\h$, such that $\phi_k = Q \psi_k$, $\forall k\in\mn$, see \cite{ba99-2}.
\end{remark}

\begin{corollary} \label{propequiv}
If $\Phi$ and $\Psi$ are equivalent frames, then $M_{(c), \Phi,\Psi}$ is invertible  and $M_{(c), \Phi,\Psi}^{-1}=M_{(\frac{1}{c}), \widetilde{\Psi},\widetilde{\Phi}}$.
\end{corollary}

Now it is natural to pose the inverse question: {\em If $M_{(c), \Phi,\Psi}^{-1}=M_{(\frac{1}{c}), \widetilde{\Psi},\widetilde{\Phi}}$, does it follow that $\Phi$ and $\Psi$ are equivalent?} We give an affirmative answer in the next theorem.

\begin{theorem}
Let $\Phi$ and $\Psi$ be frames for $\h$. The following statements are equivalent.
\begin{itemize}
\item[{\rm (a)}] $M_{(c), \Phi,\Psi}$ is invertible and $M_{(c), \Phi,\Psi}^{-1}=M_{(\frac{1}{c}), \widetilde{\Psi}, \widetilde{\Phi}}$.
\item[{\rm (b)}] $\Phi$ and $\Psi$ are equivalent frames.
\item[{\rm (c1)}] $M_{(c), \Phi,\Psi}$ is invertible and  the unique frame $\Psi^\dagger$ in \mbox{Theorem \ref{ff}} is $\widetilde{\Psi}$.
\item[{\rm (c2)}] $M_{(c), \Phi,\Psi}$ is invertible and the unique frame $\Phi^\dagger$ in \mbox{Theorem \ref{ff}} is $\widetilde{\Phi}$.
\item[{\rm (d1)}] $M_{(c), \Phi,\Psi}$ is invertible and $M_{(c), \Phi,\Psi}^{-1}=M_{(\frac{1}{c}), \widetilde{\Psi}, \Phi^d}$ for all dual frames $\Phi^d$ of $\Phi$. 
\item[{\rm (d2)}] $M_{(c), \Phi,\Psi}$ is invertible and $M_{(c), \Phi,\Psi}^{-1}=M_{(\frac{1}{c}),  \Psi^d, \widetilde{\Phi}}$ for all dual frames $\Psi^d$ of $\Psi$. 
\end{itemize}
\end{theorem}
\bp Without loss of generality, we may consider $c=1$. For a closed subspace $U$ of $\h$, the orthogonal projection on $U$ will be denoted by $P_U$.

(a) $\Rightarrow$ (b): By (i), we have $T_{\widetilde{\Psi}} U_{\widetilde{\Phi}} T_\Phi U_\Psi = Id_\h$ and hence, 
$U_\Psi T_{\widetilde{\Psi}} U_{\widetilde{\Phi}} T_\Phi U_\Psi T_{\widetilde{\Psi}} = U_\Psi T_{\widetilde{\Psi}} $.
Then $P_{\range(U_\Psi)} P_{\range(U_{\widetilde{\Phi}})}  P_{\range(U_\Psi)} = P_{\range(U_\Psi)}$, which implies that 
$ \range(U_\Psi) \subseteq \range(U_{\widetilde{\Phi}})$.

 In an analog way,  it follows that  
$  \range(U_{\widetilde{\Phi}}) \subseteq \range(U_\Psi)$. 

Therefore, $  \range(U_{\widetilde{\Phi}}) = \range(U_\Psi)$. 
This implies that $\Phi$ and $\Psi$ are equivalent.

(b) $\Rightarrow$ (c1) and (c2): Since $\psi^\dagger_n = M^{-1}(\phi_n)$, $n\in\mn$, it follows that $\Psi^\dagger$ is equivalent to $\Psi$. Therefore,  $\Psi^\dagger = \widetilde{\Psi}$. The validity of (c2) follows in a similar way.

 (c1) $\Rightarrow$ (d1) and (c2) $\Rightarrow$ (d2): Use Theorem \ref{ff}.

 (d1) $\Rightarrow$ (a) and (d2) $\Rightarrow$ (a):  clear. 
\ep

\vspace{.1in} 
For the more general case of semi-normalized symbols, it is not difficult 
to prove the following sufficient condition for validity of Equation (\ref{invmrb}). 

\begin{proposition}\label{propq2}
Let $\Phi$ and $\Psi$ be frames for $\h$, and let 
the symbol $m$ be semi-normalized. 
Assume that $M_{m, \Phi, \Psi}$ is invertible. 
If $\Psi$ is equivalent to $(m_n\phi_n)$ or 
  $\Phi$ is equivalent to $(\overline{m}_n\psi_n)$, then $M_{m, \Phi, \Psi}^{-1}=M_{1/m,\widetilde{\Psi},\widetilde{\Phi}}$.
\end{proposition}

\section{Gabor Multipliers}  \label{gabmult} 

In this section we are interested in  invertible Gabor frame multipliers, whose inverses can be written as Gabor frame multipliers, not just by frame multipliers.

\begin{theorem} \label{gab2} Let $g \in \LtRd$ and  let $(\pi(\lambda)g)_{\lambda\in\Lambda}$ be a Gabor frame for $\LtRd$. 
Let $V:\LtRd\to\LtRd$ be a bounded bijective operator. Then the following statements are equivalent. 
\begin{itemize}
\item[{\rm (${\mathcal A}_1$)}] For every $\lambda\in\Lambda$, $V \pi(\lambda) g =\pi(\lambda)Vg$.
\item[{\rm (${\mathcal A}_2$)}] For every $\lambda\in\Lambda$ and every  $f\in \LtRd$, $V\pi(\lambda)f=\pi(\lambda)Vf$ (i.e.,  $V$ commutes with $\pi(\lambda)$ for every $\lambda\in\Lambda$).
\item[{\rm (${\mathcal A}_3$)}]
$V$ can be written as a Gabor frame multiplier with the constant symbol  $(1)$. 
\item[{\rm (${\mathcal A}_4$)}]
$V^{-1}$  can be written as a Gabor frame multiplier with the constant symbol  $(1)$.
\end{itemize}
\end{theorem}
\bp
(${\mathcal A}_3$) $\Longrightarrow$ (${\mathcal A}_2$) Let $V$ be the Gabor frame multiplier $M_{(1),(\pi(\lambda)v )_{\lambda\in\Lambda},(\pi(\lambda) u)_{\lambda\in\Lambda}}$ for some $u,v\in\LtRd$. 
For every $f\in\LtRd$, every $\lambda=(\tau,\omega)\in\Lambda$ and $\lambda'=(\tau',\omega')\in\Lambda$  we have
\begin{eqnarray*}
V   \pi(\lambda) f 
 & = & \sum \limits_{\lambda'\in\Lambda} \left< \pi(\lambda) f , \pi(\lambda') u \right> \pi(\lambda')  v\\
& = & \sum \limits_{\lambda'\in\Lambda} \left<  f , e^{2\pi i \tau(\omega'-\omega)} \pi(\lambda'-\lambda) u \right> \pi(\lambda')  v\\
& = & \sum \limits_{\lambda''\in\Lambda} \left<  f , e^{2\pi i \tau\omega''} \pi(\lambda'') u \right> \pi(\lambda''+\lambda)  v\\
& = & \sum \limits_{\lambda''\in\Lambda} \left<  f , e^{2\pi i \tau\omega''} \pi(\lambda'') u \right> e^{2\pi i \tau\omega''} \pi(\lambda)\pi(\lambda'')  v \\
 &=&  \pi(\lambda)  Vf.
\end{eqnarray*} 
This statement extends the result
that the Gabor frame operator commutes with $\pi(\lambda)$; the proof uses similar techniques as in \cite[Lemma 9.3.1]{ole1}.

\v
(${\mathcal A}_2$) $\Longrightarrow$ (${\mathcal A}_1$) is obvious.

\v
(${\mathcal A}_1$) $\Longrightarrow$ (${\mathcal A}_3$) 
For every $f\in\LtRd$,
\begin{equation*} \label{xxxx}
V f = V \left( \sum \limits_{\lambda\in\Lambda} 
 \left< f , {\widetilde g}_{\lambda}\right> g_\lambda \right) =
\sum \limits_{\lambda\in\Lambda} \left< f , {\widetilde g}_{\lambda}\right>  \pi(\lambda) Vg,
\end{equation*} 
which means that $V$ can be written as a Gabor frame multiplier with symbol $(1)$.

\v
(${\mathcal A}_1$) $\Longrightarrow$ (${\mathcal A}_4$) For $\lambda\in\Lambda$, denote $h_{\lambda}=\pi(\lambda) Vg$. By what is already proved, 
$V$ can be written as the multiplier
$M_{(1),(h_{\lambda}),({\widetilde g}_{\lambda})}$. Since $(h_{\lambda})$ and $({\widetilde g}_{\lambda})$ are equivalent frames, 
Corollary \ref{propequiv} implies that 
$M^{-1}_{(1),(h_{\lambda}),({\widetilde g}_{\lambda})}=M_{(1),(g_{\lambda}),({\widetilde h}_{\lambda})}$.

\v
(${\mathcal A}_4$) $\Longrightarrow$ (${\mathcal A}_2$)
Having in mind the implication (${\mathcal A}_3$) $\Longrightarrow$ (${\mathcal A}_2$) applied to $V^{-1}$, it follows that 
$V^{-1}$ commutes with $\pi(\lambda)$,  $\forall \lambda\in\Lambda$. 
Therefore, 
$V$ also commutes with $\pi(\lambda)$, $\forall \lambda\in\Lambda$.
\ep

\begin{remark}
This result gives a nice representation and criterion for TF-lattice invariant operators \cite{feikoz1}, which correspond to condition (${\mathcal A}_2$).
Motivated by \cite{DL98}, the condition (${\mathcal A}_1$) can be considered to define 'locally TF-lattice invariant' operators. We have shown that this local property already implies the global one.   
\end{remark}

As a consequence of Theorem \ref{gab2}, the inverse of every invertible Gabor frame multiplier with constant symbol  can be written as a Gabor frame multiplier:

 \begin{corollary}
Let $V:\LtRd\to\LtRd$ be an invertible Gabor frame multiplier $M_{(1),(\pi(\lambda)v )_{\lambda\in\Lambda},(\pi(\lambda) u)_{\lambda\in\Lambda}}$. Let  $(g_{\lambda})_{\lambda\in\Lambda}=(\pi(\lambda)g)_{\lambda\in\Lambda}$ be any Gabor frame for $\LtRd$. Then $V^{-1}$ can be written as the Gabor frame multiplier  
$M_{(1),(g_{\lambda}),({\widetilde h}_{\lambda})}$, where $h_{\lambda}=\pi(\lambda) Vg$, $\lambda\in\Lambda$.
\end{corollary}

Concerning Theorem \ref{gab2}, note that if weaker assumptions on $V$ are made, then a similar  proof can be used to show the following statements.

\begin{lemma} \label{lemg} As in Theorem \ref{gab2}, let $g \in \LtRd$ and  let $(\pi(\lambda)g)_{\lambda\in\Lambda}$ be a Gabor frame for $\LtRd$. 
Let $V:\LtRd\to\LtRd$ be an operator. Consider the condition
\begin{itemize}
\item[{\rm (${\mathcal A}_3'$)}] $V$ can be written as a Gabor Bessel multiplier with a constant symbol  $(1)$.  
\end{itemize} Then the following statements hold.
\begin{itemize}
\item[{\rm (i)}] If $V$ is bounded, then  {\rm (${\mathcal A}_3$)  $\Rightarrow$ (${\mathcal A}_3'$)  $\Leftrightarrow$   (${\mathcal A}_2$) $\Leftrightarrow$   (${\mathcal A}_1$)};

\item[{\rm (ii)}] If $V$ is  bounded and surjective, then {\rm (${\mathcal A}_3$)  $\Leftrightarrow$  (${\mathcal A}_3'$)
  $\Leftrightarrow$ (${\mathcal A}_2$) $\Leftrightarrow$   (${\mathcal A}_1$)}.
\end{itemize}
\end{lemma}

\sloppy
So, when a Gabor Bessel multiplier $M_{(1),(\pi(\lambda)v )_{\lambda\in\Lambda},(\pi(\lambda) u)_{\lambda\in\Lambda}}$ is bounded and surjective, it can always be written as a Gabor frame multiplier for appropriate frames.
Note that if $V$ is the Gabor Bessel multiplier $M_{(1),(\pi(\lambda)v )_{\lambda\in\Lambda},(\pi(\lambda) u)_{\lambda\in\Lambda}}$ for some $u,v\in\LtRd$ and $V$ is a bounded surjective (resp. invertible) operator, then $(\pi(\lambda) v)_{\lambda\in\Lambda}$ is a frame (resp. $(\pi(\lambda) u)_{\lambda\in\Lambda}$ and $(\pi(\lambda) v)_{\lambda\in\Lambda}$ are frames) for $\LtRd$ (see e.g. \cite{bas} (resp. see \cite{balsto09new})). 

\section{Numerical Visualization of Results}
Let us come back to the example in Figure \ref{fig:examp1}. 
Here we use the same signal $f$, Gabor frame $\Psi$, and symbol $m$, as in Figure \ref{fig:examp1}. 
Note that all the elements of the symbol $m$ fulfill  $m_{n,k}\in\{1,10^{-2}\}$ and denote  $M = M_{m, \widetilde \Psi, \Psi}$.
Since $m$ is semi-normalized, the multiplier $M$ is analytically invertible \cite[Prop. 4.3]{balsto09new}. 
However, the operator is badly conditioned, the condition number is around $99$.
As mentioned before, the signal $f$ is approximately $2$ seconds long, using a sampling rate of $44100$. So the signal is a $128148$-dimensional vector. 

Starting from $g = M f$, we compare two approaches numerically, 
\begin{enumerate}
\item a 'naive' inversion $ \dot{f} = M_{1/m, \widetilde \Psi, \Psi} g$ (corresponding to the approach raised in [Q2]),
\item and the 'iterative' inversion  $\ddot{f} = M^{-1} g$. 
For numerical efficiency and in particular, for memory constraints, we use the iterative inversion 
in LTFAT, using a conjugate gradient method (for $M^* M$). 
 Note that by Theorem \ref{ff}, 
$M^{-1} g$ corresponds to $M_{1/m, \Psi^\dagger, \Psi} g$, 
 where $\psi_{n,k}^\dagger = \left( M^{-1} \left(m_{n,k} \tilde \psi_{n,k} \right) \right)$.
\end{enumerate}  

\begin{figure}[ht]
\includegraphics[width=1\textwidth]{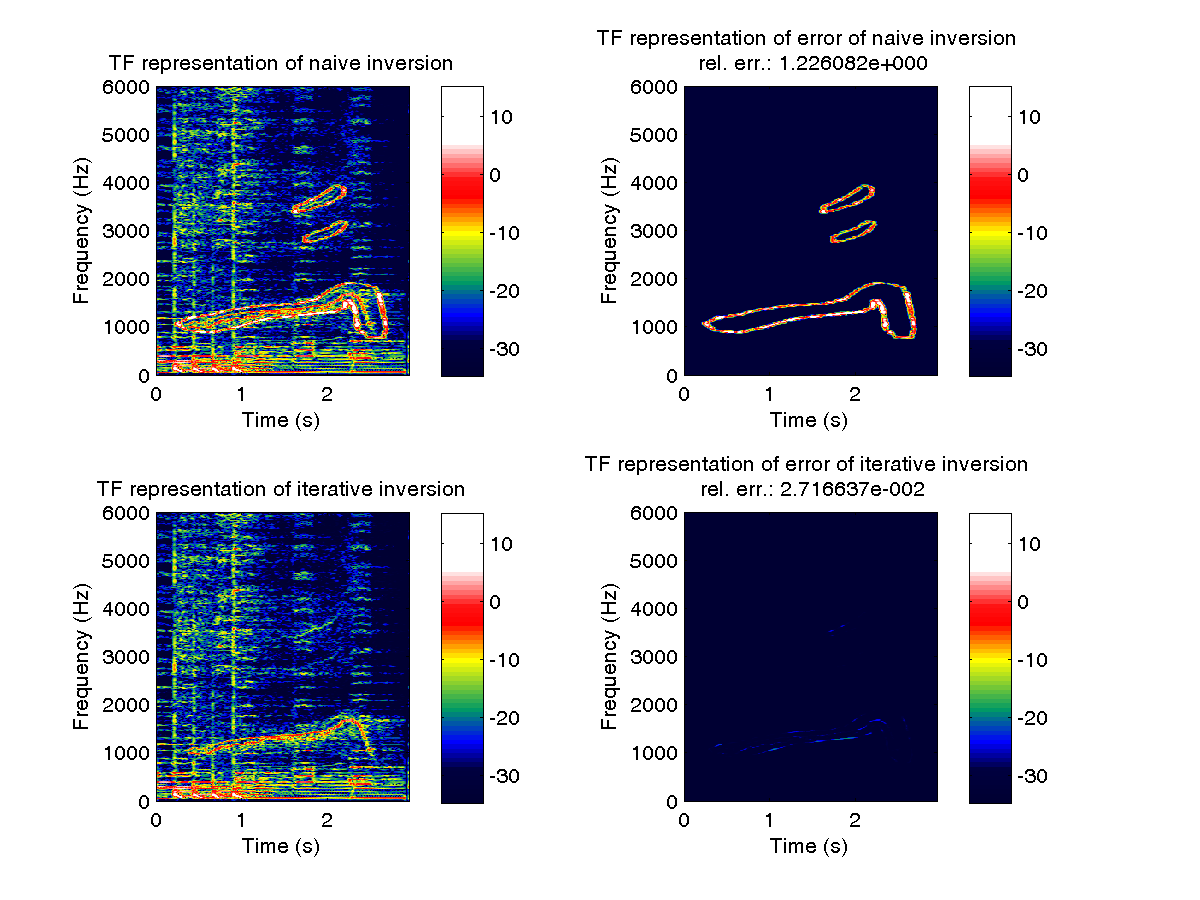} 
\caption{\label{fig:vis2} {\em Inversion of multipliers.} 
 Time-Frequency representation of (TOP LEFT) the result of the 'naive' inversion $\dot{f}$ (TOP RIGHT) the error of the 'naive' inversion, i.e.  $\dot{f} - f$ (BOTTOM LEFT) the iterative inversion $\ddot{f}$ (BOTTOM RIGHT) the error of the iterative inversion  $\ddot{f} - f$.} 
\end{figure}

For results see Figure \ref{fig:vis2}. 
Clearly the naive approach has strong artifacts. 
The error is especially big at the boundaries of the constant region of the symbols. 
The chosen atoms are well localized in time-frequency, so that within the interior of the constant regions, this inversion works well. 
This could be expected as we have shown in Corollary \ref{propequiv} that constant symbols allow this kind of inversion for equivalent frames (so, in particular for $\Psi$ and $\widetilde \Psi$).

The iterative inversion worked well with an error of $3$ \%. This could, naturally, be decreased by investing more calculation time. But also in the chosen setting for the iterative inversion ($100$ iterations in {\sf iframemul} \cite{soendxxl10}) no difference can be seen in the time-frequency representation, as well as no audible difference can be detected.

Similar results can also be shown for other redundancies and other sound files.

\v
 {\bf Acknowledgments}
 
The work on this paper was supported by the WWTF project MULAC ('Frame Multipliers: Theory and Application in Acoustics; MA07-025)  and 
the Austrian Science Fund (FWF) START-project FLAME ('Frames and Linear Operators for Acoustical Modeling and Parameter Estimation'; Y 551-N13). 
  The authors are thankful to H. Feichtinger and D. Bayer for their valuable comments, as well as P. Soendergaard and Z. Prusa for help with the implementations.
   They also thank NuHAG (University of Vienna) for the possibility to use their online database (http://www.nuhag.eu). 
The second author is grateful for the hospitality of the Acoustics Research Institute and the support from the projects MULAC and FLAME.

\providecommand{\MR}{\relax\ifhmode\unskip\space\fi MR }

\providecommand{\MRhref}[2]{%
  \href{http://www.ams.org/mathscinet-getitem?mr=#1}{#2}
}
\providecommand{\href}[2]{#2}

\end{document}